\input amstex
\documentstyle{amsppt}

\loadbold
\topmatter
\pageheight{7.5in}

\title  The Fine Structure of the Kasparov Groups III: \\
Relative Quasidiagonality
      \endtitle
 \rightheadtext{ Fine Structure III }
\author    C. L. Schochet         \endauthor
 \affil
 Mathematics Department   \\
 Wayne State University    \\
 Detroit, MI 48202      \\ \\
 \endaffil
\address       Wayne State University     \endaddress
\email     {    claude\@math.wayne.edu }     \endemail
 \date {July 10, 2001    }          \enddate
\endkeywords
  \endsubjclass
\abstract { Quasidiagonality was introduced by P.R. Halmos for operators and quickly
generalized to $C^*$-algebras. D. Voiculescu   asked
how quasidiagonality (in its various forms) was related to topological
invariants.
N. Salinas     systematically studied the  topology of the Kasparov
groups
$KK_*(A,B) $ and showed that this topology is
related to relative quasidiagonality.
In this paper we
 identify $QD(A,B)$, the   quasidiagonal classes
in $  KK_1(A,B)  $,
in terms of $K_*(A)$ and $K_*(B)$,
 and we use these results in various applications.
Here is our central result.
Let $\widetilde{\Cal N} $ denote the category of separable 
nuclear $C^*$-algebras 
which satisfy the Universal Coefficient Theorem.
\medbreak
\proclaim {Theorem  } Suppose that $A \in \widetilde{\Cal N} $ and
$A$ is quasidiagonal relative to $B$.
Then there is a natural isomorphism
$$
QD(A,B) \,\cong\,   Pext_{\Bbb Z}^1( {K_*(A)},{K_*(B)}) _{0}     .
$$
\endproclaim
\medbreak
Thus for $A \in \widetilde{\Cal N} $   quasidiagonality of $KK$-classes is indeed a
topological invariant.
We also settle a question raised by L.G. Brown on the relation
between relative quasidiagonality and the kernel of the natural map
$$
\theta ^*  : Ext_{\Bbb Z}^1 ( K_*(A) , K_*(B) )
\longrightarrow Ext_{\Bbb Z}^1 (   K_*(A)_t , K_*(B)  ) .
$$
where $K_*(A)_t $ denotes the torsion subgroup of $K_*(A)$.
Finally we establish a converse to a theorem of Davidson, Herrero, and
Salinas, giving conditions under which the quasidiagonality
of $A/\Cal K$ implies the quasidiagonality of the associated 
representation of $A$.
}\endabstract

\endtopmatter

\document
\magnification = 1200
\pageno=1

\def\tensor{\otimes}

\def\KK #1.#2 { KK_*( #1,#2) }
 \def\KKgraded #1.#2.#3 { KK_{#1}( #2, #3) }

\def\ext #1.#2 {Ext _{\Bbb Z}^1( {#1} , {#2} ) }
\def\pext #1.#2 {Pext _{\Bbb Z}^1( {#1} , {#2} ) }
\def\usualpext {\pext{K_*(A)}.{K_*(B)} }
\def\usualext {\ext{K_*(A)}.{K_*(B)} }

\def\hom #1.#2 {Hom_{\Bbb Z} (#1 ,  #2  ) }

\def\starone{{*-1}}
\def\usualhom {\hom{K_*(A)}.{K_*(B)} }

 \newpage

\newpage\beginsection {1. Introduction: Quasidiagonality and $KK$-theory}

This is the third  in a series   of papers in
which the topological structure of the Kasparov groups, systematically
studied first  by Salinas, is developed and put to use. The first two
papers \cite {20, 21} are  devoted to general structural results and
serve as the theoretical background for the present work,
which centers about quasidiagonality. From the point of view
of \cite {20, 21}, this is an exploration of the closure of zero
in the Kasparov groups, which we have termed the {\it {fine
structure}} subgroup.

Quasidiagonality was defined by P.R. Halmos \cite {12} in 1970. A bounded
operator on    Hilbert space is {\it {quasidiagonal }} if it is a compact
perturbation of a block-diagonal operator. This soon was generalized to
$C^*$-algebras. Quasidiagonality is
thus a finite dimensional approximation property.  It is  not
well understood.

L.G. Brown, R.G. Douglas, and P.A. Fillmore \cite {6} first recognized
that the study of
quasidiagonality for operators and for $C^*$-algebras might
be approached by topological methods.  They topologized
their functor $\Cal Ext(X)$ (which is known now to be
isomorphic to the Kasparov group
$\KKgraded 1.{C(X)}.{\Bbb C} $)
and announced that the closure of zero corresponded to the quasidiagonal
extensions. L.G. Brown pursued this
theme, particularly in \cite {4} (cf. \S\S 6,7).

Salinas \cite {15}   studied the topology on the Kasparov
groups $\KKgraded *.A.B $ and showed that this topology is
related to   relative quasidiagonality.
The quasidiagonal classes $QD(A,B)$  (defined precisely in
\S 2) constitute a certain   subgroup
of $KK$-theory:
$$
QD(A,B) \,\subseteq\, \KKgraded 1.A.B    .
$$
 If $A$ is in the category $\widetilde{\Cal N} $ of 
 separable nuclear\footnote{We require nuclearity for two reasons.
 First, we need at least $K$-nuclearity so that the Kasparov
 product will be continuous (cf. 1.2). Second, we apparently 
 need nuclearity in order to satisfy the hypotheses of Salinas's 
 result identifying the quasidiagonal elements with the closure
 of zero in the Kasparov groups (cf. 2.2). It seems possible
 that if one restricts attention to extensions that have additive
 inverses so that the identification 
$$
\Cal Ext_*(A,B) \,\cong\, \KKgraded *.A.B       
$$
 holds, 
 then his result might generalize to the $K$-nuclear 
 setting.} 
   $C^*$-algebras which 
 satisfy the Universal Coefficient Theorem
 \footnote {To be precise, $\widetilde {\Cal N}$ is the full subcategory
 of the category of separable nuclear 
  $C^*$-algebras $A$  that have the property that
   the UCT is satisfied
 for every pair $(A,B)$, where $B$ ranges over all separable
 $C^*$-algebras $B$.  Equivalently, we could require that 
 $A$ is $KK$-equivalent to a commutative $C^*$-algebra. The
  category $\widetilde{\Cal N} $ includes the bootstrap category $\Cal N$,
  by the UCT of Rosenberg and the author \cite {14}.  
 The bootstrap
category   $\Cal N$ \cite {16, 14}
is the smallest full subcategory of separable nuclear
$C^*$-algebras which contains all separable Type I $C^*$-algebras
and which is closed under strong Morita equivalence, direct limits, extensions,
and crossed products by $\Bbb Z$ and $\Bbb R$. We may also require that if
$J$ is an ideal of $A$ and if $J$, $A \in \Cal N$ then so is $A/J$, and if
$A$ and $A/J$ are in $\Cal N$ then so is $J$. The hope is 
that every separable nuclear $C^*$-algebra is $KK$-equivalent to 
a commutative $C^*$-algebra. }
    then  more can be said.
   The UCT is a
natural short exact sequence
$$
0 \rightarrow \ext {K_*(A )}.{K_*(B) }
\overset\delta\to\longrightarrow \KK {A }.{B}
\overset\gamma\to\longrightarrow   \hom {K_*(A )}.{K_*(B) }
\rightarrow 0
$$
which splits unnaturally and thus computes $\KKgraded *.A.B $ in terms
of $K_*(A)$ and $ K_*(B)$.   In particular, it identifies a canonical
subgroup of $\KKgraded *.A.B  $, namely
$$
\ext {K_*(A )}.{K_*(B) }  \,\overset\delta\to\cong\, Ker(\gamma)
\hookrightarrow   \KKgraded *.A.B
$$
Henceforth we generally suppress mention of the map
$\delta $.\footnote {   Note, though, that the map
$\delta $ has degree one and so the elements of degree  zero in
the group $       \ext {K_*(A )}.{K_*(B) }      $, denoted by
$  {  \ext {      K_*(A)    }.{     K_*(B)    } }_0$,
 are contained in  $\KKgraded 1.A.B $.          }

Salinas has shown \cite {15,  5.1}   that
if $A$ is quasidiagonal relative to $B$ then 
$$
QD(A,B) \subseteq   Ker(\gamma )
$$
and in fact   \cite {15, Theorem 5.2a
 as reformulated   by M. Dadarlat (private
communication)}  that if $A \in \widetilde{\Cal N} $ so that the UCT holds
and identifies 
$$
Ker(\gamma ) \,\cong\, \ext {K_*(A )}.{K_*(B) } 
$$
then
$$
QD(A,B) \subseteq   \pext {K_*(A )}.{K_*(B) } _0
$$
where $\pext G.H $ is the subgroup of $\ext G.H $ consisting of
pure extensions.\footnote {A subgroup $H \subseteq J$
is said to be {\it pure} if for each $n \in \Bbb N$,
$$
nH = H \cap nJ  
$$
and an extension of abelian groups
$$
0 \to H  \to J  \to G  \to 0
$$
is said to be {\it pure } if $H$ is a pure subgroup of $J$. If $H$ is a
direct summand of $J$ then $H$ is a pure subgroup,
 but the  most interesting cases
involve non-split pure extensions. For example, $tJ$, the torsion
subgroup of $J$,  is always a pure
subgroup of $J$   but it is not necessarily a
direct summand of $J$.
See \cite {11,  \S 53}  
  and \cite{22}.}
\medbreak

In the first two papers  in this series \cite {20, 21} we demonstrated
the following facts:

\medbreak

\flushpar{\bf 1.1}.\quad There is a natural structure of a pseudopolonais\footnote{
A topological space is {\it {polonais}} if it is separable, complete,
and metric. If it is a topological group then we insist that the
metric be invariant. A {\it {pseudo\-polonais}} group is a separable
  topological group with invariant 
pseudo\-metric whose Hausdorff quotient group is polonais.}
group
on $\KKgraded *.A.B  $ \cite {20, 6.2}.
\medbreak
\flushpar{\bf 1.2}\quad The Kasparov pairing is jointly continuous with respect
to this topology, provided that all $C^*$-algebras
 which appear in the first variable are $K$-nuclear \cite {20, 6.8}.
\medbreak
 
\flushpar{\bf 1.3}.\quad
 If $K_*(A)$ is finitely generated
then $\KKgraded *.A.B $ is polonais. \cite {20, 6.2}
\medbreak

\flushpar{\bf 1.4}.\quad The index map
$$
\gamma : \KKgraded *.A.B  \to \hom {K_*(A)}.{K_*(B)}
$$
is continuous.
If $Im(\gamma ) $ is closed (e.g., if $\gamma $ is onto), then
$\gamma $ is an open map.
 If $\gamma $ is an algebraic isomorphism then it is an
isomorphism of topological groups.  \cite {20, 7.4}
\medbreak

\flushpar{\bf 1.5}.\quad The Universal Coefficient Theorem short exact sequence
is a sequence of pseudo\-polonais groups and  each of the
splittings of the UCT constructed in \cite {14} is a topological splitting.
\cite {21, 4.5}

\medbreak\flushpar{\bf 1.6}.\quad If  $A \in \widetilde{\Cal N} $
  then there is a natural isomorphism
$$
Z_*(A,B) \cong \pext {K_*(A)}.{K_*(B)}
$$
where $Z_*(A,B)$ denotes the closure of zero in the
group $\KKgraded *.A.B $. \cite {21, 3.3}
\medbreak

\flushpar{\bf 1.7}.\quad If $A \in \widetilde{\Cal N} $ then it has an associated
$KK$-filtration diagram for $(A,B)$ which is functorial
into the category of pseudopolonais groups \cite{21}.\footnote{
We proved 1.6 and 1.7 in \cite{21} for $A \in \Cal N$. However, 
the results hold for $A \in \widetilde{\Cal N} $ as may easily be seen by inspection.}
In particular, the Milnor and Jensen sequences take values
in this category and are natural with respect to both
$A$ and $B$. 
\medbreak

The following theorem is the  most important result in this paper; all of
our applications flow from it. Salinas \cite{15}
  proved the theorem under certain stringent
assumptions on $A$ and $K_*(A)$ which we have removed.
\medbreak

\proclaim {Theorem 2.3} Suppose that $A \in \widetilde{\Cal N} $ and
$A$ is quasidiagonal relative to $B$.
Then there is a
natural isomorphism  of topological groups
$$
   QD(A,B)  \,\cong\,  \pext  {K_*(A)}.{K_*(B)} _0
$$
  regarded as topological subgroups
of $\KKgraded *.A.B   $ via
the canonical inclusion $\delta $ in the UCT.
\endproclaim
\medbreak

Note   as an immediate consequence of this theorem that  
quasidiagonality 
of extensions
is a topological invariant for $A \in \widetilde{\Cal N} $, answering
 the relative form of
a question of D. Voiculescu \cite{24,\, 25}. For instance (see Theorem 3.5), if
$$
\CD
0 @>>> B\tensor\Cal K  @>>>  E_\tau   @>>>  A  @>>>  0 \\
@. @VV1V    @VV{e_\tau }V   @VV{\tau }V  \\
0 @>>> B\tensor\Cal K  @>>> \Cal L(H_B)  @>>> \Cal Q(H_B) @>>> 0
\endCD   
$$
is an essential extension classified by $\tau $ such that $A \in \widetilde{\Cal N} $
  and $A$ is quasidiagonal
 relative to $B$, then $e_\tau (E_\tau )$ 
 is $B$-quasidiagonal  
if and only if both
of the following topological conditions hold:
\roster
\item $\gamma (\tau )  = 0: K_*(A) \to K_*(B)$, and
\medbreak
\item
 $$
\tau \in \underset{n}\to\bigcap \,\, n\ext {K_*(A)}.{K_*(B)}   _0
= \pext {K_*(A)}.{K_*(B)}   _0 .
$$
\endroster
 If $K_*(A)$ is torsionfree then condition 2) is automatically satisfied,
so  $e_\tau (E_\tau )$ is $B$-quasidiagonal  
 if and only if $\gamma (\tau ) = 0$.
Thus when $A \in \widetilde{\Cal N} $ with $K_*(A)$ torsionfree, the index
invariant $\gamma $ is a complete obstruction to relative quasidiagonality.

\medbreak

The remainder of the paper is organized as follows.
\medbreak

\medbreak
In Section 2  the definitions of quasidiagonality are recalled and
Theorem 2.3 is established.

              Theorem 2.3 has several corollaries
which are developed in \S 3. Here
is one. Suppose that  $f: A \to B$ so that 
$$
[f] \in \KKgraded 0.{A}.{B} \overset{\beta ^A}\to\longrightarrow
 \KKgraded 1.{SA}.B 
$$
where $\beta ^A$ denotes the Bott periodicity isomorphism in the 
first variable.
 When
is $\beta ^A([f])$  a quasidiagonal class? It is
 easy to show that the following are
necessary conditions:
\roster
\item The induced homomorphism
$
f_* : K_*(A) \to K_*(B) $  is trivial; and,
\medbreak
\item
The associated short exact sequence
$$
0 \to  K_*(SB)  \to K_*(Cf)  \to  K_*(A)  \to  0
$$
is pure exact, where $Cf$ denotes the mapping cone of $f$.
\endroster
\medbreak
\flushpar It is shown  that
 (for $A \in \widetilde{\Cal N} $) these conditions are also sufficient.
If $K_*(A)$ is torsionfree then (2) is automatic, and so $\beta ^A([f])$
is a quasidiagonal class if and only if $f_* = 0$.
 The section
concludes with Theorem 3.5, which demonstrates the topological nature
of relative quasidiagonality in a concrete manner.

\medbreak The remaining four sections are devoted to applications.
\medbreak

Section 4 is devoted to the application of some of the theory of
infinite abelian groups to obtain results on quasidiagonality.
 \medbreak

In Section 5 we answer a
question raised by L.G. Brown \cite {4} concerning the relation
between quasidiagonality and the kernel of the map
$$
\theta ^*  : \ext {K_*(A)}.{K_*(B)} \longrightarrow \ext {K_*(A)_t}.{K_*(B)}  .
$$
induced by the natural inclusion of $K_*(A)_t$, the torsion subgroup
of $A$, into $K_*(A)$.
\medbreak

Section  6  deals with another result of L.G. Brown.
 Brown constructed \cite {4}
an operator $T$ which was not quasidiagonal but such that $T\oplus T$ was
quasidiagonal. In Section 7 we analyze such phenomena.
\medbreak

Section 7 presents    a converse to a theorem of Davidson, Herrero,
and Salinas  which deals with conditions under which the quasidiagonality
of $A/\Cal K$ implies the quasidiagonality of $A$.

\medbreak
It is a pleasure to acknowledge helpful correspondence and conversations
regarding quasidiagonality and
abelian groups with L.G. Brown, N. Brown, M. Dadarlat,     H. Lin,
T. Loring,   N. Salinas, and D. Voiculescu,
with a special thanks to N. Salinas for his help and encouragement.

\medbreak

In this paper  all $C^*$-algebras
are assumed separable  
 with the obvious exceptions of
multiplier algebras $\Cal M (B\tensor\Cal K) $ and their quotients.
 Whenever we speak of quasidiagonal classes in $\KKgraded 1.A.B $
it is understood that 
$B$ is separable,
$B\tensor\Cal K$ has a countable approximate
unit consisting of projections,\footnote {Dadarlat calls this 
property 
{\it {stably unital}}. Note that if $B$ is unital and $\{p_i\} 
\subset \Cal K$ is a countable approximate unit for $\Cal K$ 
consisting of projections, then $\{1 \tensor p_i\}$ is
  a countable approximate unit for $B\tensor \Cal K$ 
consisting of projections.
Thus if $B$ is unital then $B \tensor\Cal K$ is stably unital.
  } and  
$A$ is quasidiagonal relative to $B$.

 An isomorphism  of topological groups is an algebraic isomorphism
which is  a homeomorphism of topological spaces.
We use the topologists' notation for graded abelian groups. For example,
$$
\usualpext _0 = \pext{K_0(A)}.{K_0(B)} \oplus \pext{K_1(A)}.{K_1(B)} .
$$
\medbreak

Note: This paper,  and to a small extent \cite {20, 21},   replace
and very substantially extend
the preliminary preprint entitled \lq\lq Continuity of the Kasparov pairing
and relative quasidiagonality\rq\rq which will not appear.


\newpage\beginsection {2. Quasidiagonality}

Our description of the various definitions of quasidiagonality
leans heavily upon the remarkable survey paper of
Nathaniel P. Brown \cite{7}
and from the paper of Marius Dadarlat \cite{9}.
 We are most grateful
to them for clarifying these issues.

Halmos \cite{12} introduced the notion of a quasidiagonal operator.
A bounded linear operator
on a separable Hilbert space
$T \in \Cal L(H)$ is a {\it {block diagonal}}
operator if there exists a countable approximate unit consisting
of projections, that is,  an increasing sequence of finite rank
projections $P_1 \leq P_2 \leq P_3 \dots $ converging to the identity
in the strong operator topology, which is central with respect to $T$:
$$
P_nT - TP_n = 0\qquad  \forall\, n.
$$
  An operator $T \in \Cal L(H)$ is
{\it {quasidiagonal}}
 if there exists a countable approximate
unit consisting of projections
$\{P_n\}$
 which is quasicentral with respect
to $T$:
$$
{\underset{n \to \infty}\to{lim}} || P_nT - TP_n || = 0.
$$
 The sum of a compact operator
 and a block diagonal operator is   quasidiagonal, and Halmos
  proved that in fact every quasidiagonal operator has this form.

The concept extends to $C^*$-algebras as follows.
Suppose that $B$ is a separable $C^*$-algebra. 
Let $H_B = B\tensor H $.
We write 
$$
\Cal L(H_B) = \Cal M(B\tensor\Cal K)\qquad {\text{and}} \qquad
 \Cal Q(H_B) = \Cal M(B\tensor\Cal K)/B\tensor\Cal K .
 $$ 
 A separable\footnote{This definition is not correct if $E$ is not
separable; see \cite{7} for the correct definition.} 
 subset $E \subset \Cal L(H_B)$
is called a {\it {$B$-quasidiagonal set}}
\footnote{Dadarlat calls this 
\lq\lq quasidiagonal\rq\rq\, but we prefer to keep track of
$B$.}
if there
exists a countable approximate unit $\{ p_n\}$ of $B\tensor\Cal K$ 
consisting of   projections
$p_1 \leq p_2 \leq p_3 \dots $
which is quasicentral with respect to each $a \in E$:
$$
{\underset{n \to \infty}\to {lim}} || p_n a - a p_n || = 0
\qquad \forall\, a \in E.
$$
A representation $\rho : A \to \Cal L(H_B) $ of a 
separable
$C^*$-algebra $A$  
is
said to be a
{\it{ $B$-quasidiagonal representation}} 
if the set $\rho (A) $ is a $B$-quasidiagonal set.
If a separable $C^*$-algebra $A$ has a faithful 
essential\footnote{That is, the induced homomorphism 
$A \overset\rho\to\longrightarrow  \Cal L(H_B) \to \Cal Q(H_B)$
is faithful, or, equivalently, $\rho (A) \cap (B\tensor\Cal K )
= \{ 0 \}$.}
$B$-quasidiagonal representation then $A$ is said to be 
a {\it {$B$-quasidiagonal}} $C^*$-algebra.

Note that a set is $\Bbb C$-quasidiagonal if and only if 
there exists a countable approximate unit consisting of projections
in $\Cal K$ such that 
each operator in the set is quasicentral with respect to the 
countable approximate unit. In other words, each operator must 
be quasidiagonal in the classical sense, and there must be 
a countable approximate unit consisting of projections which 
works for every operator.
 We   write {\it quasidiagonal} 
rather than {\it $\Bbb C$-quasidiagonal}.

Salinas \cite{15, 4.3} shows that if $A$ is quasidiagonal
(say via a quasidiagonal representation $\rho $) and
$B$ has a countable approximate unit consisting of projections
then $A$ is $B$-quasidiagonal, since we may easily construct
the requisite countable approximate unit in $B\tensor \Cal K$
which is quasicentral with respect to the representation
$$
A \overset\rho\to\longrightarrow \Cal L(H) \cong \Cal M(\Cal K) 
\longrightarrow \Cal M(B \tensor \Cal K) \cong \Cal L(H_B) .
$$

Every commutative $C^*$-algebra is quasidiagonal- this is
a consequence of the spectral theorem. It is easy to show
that any AF algebra is also quasidiagonal. As subalgebras
of quasidiagonal algebras are obviously quasidiagonal, it
follows that any $C^*$-algebra which embeds in an AF
algebra is itself quasidiagonal. For example, this implies
that the irrational rotation $C^*$-algebras are quasidiagonal.

The unitalization of a quasidiagonal $C*$-algebra is quasidiagonal,
as is the product and minimal tensor product of quasidiagonal
$C^*$-algebras. Quasidiagonality does {\it not} pass to
quotients in general. We return to this point in \S 7.

Any quasidiagonal Fredholm operator must have Fredholm index
zero. Halmos used this fact to show that the unilateral
shift is not quasidiagonal. Then any $C^*$-algebra containing
the shift   cannot
 be quasidiagonal. More generally, quasidiagonal $C^*$-algebras
must be stably finite - they and matrix rings over them may
not contain proper isometries.

Voiculescu proved \cite{25} that if $A$ and $B$
are homotopy equivalent $C^*$-algebras then  if one of them
is quasidiagonal then the other must be as well. 
Thus for instance $CA$,
the cone on any $C^*$-algebra, is quasidiagonal, being homotopy
equivalent to $0$, and the suspension $SA$ is quasidiagonal,
 as it is a subalgebra of $CA$. Salinas extended this to the
 $B$-quasidiagonality setting. We state his result formally as part
 of Proposition 2.1.

Blackadar and Kirchberg \cite{3, 4}  introduce 
the class of NF algebras and demonstrate that this class
coincides with the class of 
separable nuclear quasidiagonal $C^*$-algebras.

Suppose that  an injection 
$$
\tau : A \longrightarrow
 \Cal Q(H_B) 
$$
classifies an extension of $C^*$-algebras. Taking pullbacks yields the 
corresponding extension
$$
0 \to B\tensor\Cal K \to E_{\tau } \to A \to 0
$$
together with the canonical faithful representation 
$$
e_{\tau } : E_{\tau } \longrightarrow \Cal L(H_B).
$$
We say that this extension is a {\it quasidiagonal extension } if 
$e_{\tau } $ is a $B$-quasidiagonal representation. 
Equivalently, 
the extension is quasidiagonal if 
there is an approximate unit consisting of projections
in $B\tensor\Cal K$ which is quasicentral in $E_\tau $. 
Salinas shows that this
 property
depends only upon the equivalence class 
$$
[\tau ] \in \KKgraded 1.A.B
$$ 
and hence it makes sense to speak of a  {\it{quasidiagonal class}} in 
$\KKgraded 1.A.B $.   

Suppose given an essential extension 
$$
0 \to B \tensor\Cal K \to E_\tau \to A \to 0
$$
together with the canonical faithful representation 
$e_\tau : E_\tau \to \Cal L(H_B) $, and suppose given 
a $*$-homomorphism $f: A' \to A$ which is an injection. Then the induced 
extension $f^*\tau $ is obtained by the pullback diagram
$$
\CD
0 @>>>  B\tensor\Cal K   @>>>   E_{f^*\tau }  @>>> A' @>>> 0  \\
@.  @VVV   @VV{\tilde f}V   @VV{f}V  \\
0 @>>>  B\tensor\Cal K   @>>>   E_{\tau }  @>>> A @>>> 0
\endCD
$$
with $\tilde f$ also an injection.
Then $e_{f^*\tau }  : E_{f^*\tau} \to \Cal L(H_B) $ 
is given by  the composition 
$$
E_{f^*\tau } \overset{\tilde f}\to\longrightarrow E_\tau 
\overset{e_\tau }\to\longrightarrow \Cal L(H_B)  .
$$
If $e_\tau $ is a $B$-quasidiagonal representation then 
it is clear from the construction that 
  the canonical faithful representation $e_{f^*\tau }$
is also  $B$-quasidiagonal.\footnote{Using Salinas's result 
(last part of Proposition 2.1) we can see that it is not necessary to stipulate 
that $f$ be an injection, once it is established that the 
homomorphism   
$$
f^* : \KKgraded *.A.B  \longrightarrow \KKgraded *.{A'}.B
$$
is continuous.}

A separable nuclear $C^*$-algebra
$A$ is said to be {\it {quasidiagonal relative to $B$}} if 
the class of the trivial extension $0 \in \KKgraded 1.A.B $
is quasidiagonal. 
If $A$ is quasidiagonal relative to $B$ then trivial extensions
are quasidiagonal, obviously,
but there may be no other quasidiagonal classes.\footnote{For example,
if $K_*(A)$ is a direct sum of cyclic groups then Theorem 2.3 implies
that every quasidiagonal extension is trivial. We discuss such matters 
in some length in later sections.}

N. Brown remarks \cite{7, 8.2} that it is possible to have an 
extension
$$
0 \to B\tensor\Cal K   \longrightarrow E \longrightarrow A \to 0
$$
with
 $E$ 
 quasidiagonal relative to $B$ without either $B$ or $A$
being quasidiagonal. However, given a separable $C^*$-algebra $D 
\subset  \Cal L(H)$
then $D$ is a quasidiagonal set if and only if 
the extension
$$
0 \to \Cal K \longrightarrow C^*\{D,\Cal K\}  \longrightarrow A \to 0
$$
is a quasidiagonal extension, where  $ C^*\{D,\Cal K\} $ denotes
the $C^*$-algebra generated by $D$ and by $\Cal K$ in $\Cal L(H)$.
 We expand upon this remark as follows.

\proclaim{Proposition 2.1} Let $A$ and $B$ be   separable
 nuclear   $C^*$-algebras.
  Then the 
 following are equivalent:
 \medbreak
 \roster 
 \item $A$ is quasidiagonal relative to $B$.
 \medbreak
 \item $A$ is $B$-quasidiagonal.
 \endroster\medbreak
  If these hold then they also hold for $*$-subalgebras of $A$.
   In addition, if $A$ is $B$-quasidiagonal and $A$ is homotopy 
  equivalent to $A'$ then $A'$ is $B$-quasidiagonal.
 \endproclaim
 
 \medbreak
 
 \demo{Proof} Suppose first that $A$ is quasidiagonal relative to $B$.
 Then there exists an essential  extension $\tau $ 
 representing $0 \in \KKgraded 1.A.B $
 and a 
 commuting classifying diagram
 $$
 \CD
 0 @>>> B\tensor\Cal K @>>> E_\tau @>>> A @>>> 0  \\
 @. @VV1V   @VVe_\tau V   @VV\tau V  \\
 0 @>>> B\tensor\Cal K @>>> \Cal L(H_B) @>>> \Cal Q(H_B)  @>>> 0  
 \endCD
 $$
 and the map $\tau $ lifts to a $*$-homomorphism 
 $\tilde\tau : A \to \Cal L(H_B)$
 since $[\tau ] = 0$. By assumption $e_\tau $
 is a $B$-quasidiagonal representation, and hence $e_\tau (E_\tau )$
 is a $B$-quasidiagonal set. The representation $\tilde\tau $ 
 has range contained in  $e_\tau (E_\tau )$ and hence
 $\tilde\tau (A)$ is a $B$-quasidiagonal set. Thus
  $\tilde\tau $
 is an essential $B$-quasidiagonal representation of $A$. So 
 $A$ is $B$-quasidiagonal.
 
In the other direction, suppose that $A$ is $B$-quasidiagonal. Then
there exists an essential  $B$-quasidiagonal representation
 $\tilde\tau : A \to \Cal L(H_B)$. Let 
 $$
 E = C^*\{ \tilde\tau (A) , B\tensor\Cal K \}
 $$
 and let $e : E \to \Cal L(H_B) $ denote the natural inclusion. 
 Then $e$ is a $B$-quasidiagonal representation and 
 $$
 0 \to B\tensor\Cal K \to E \to A \to 0
 $$
 is an essential quasidiagonal extension which is split by 
 the map $\tilde\tau $. Thus $A$ is quasidiagonal relative to $B$.
 
 If the properties hold then one uses the fact pointed
 out previously that the 
 restriction of a $B$-quasi\-diagonal extension to a subalgebra 
 of $A$ is again a $B$-quasi\-diagonal extension.
 The final statement was established by 
 Salinas \cite{15, Theorem 4.8}. \qed
 \enddemo

\medbreak

Let $QD(A,B) $ denote
 the set of quasidiagonal  classes of $\KKgraded 1.A.B     $.
This set is non-empty if and only if $A$ is quasidiagonal
relative to $B$ by \cite{15} Theorem 4.4.   
Salinas has shown that the quasidiagonal elements may be described in terms of
the topology on $\KKgraded *.A.B $.  For $B = \Cal K$ the following theorem
was established by L.G. Brown \cite {4, p. 63, Remark 1.}
\medbreak

\proclaim {Theorem 2.2 (Salinas   \cite {15, Theorem 4.4})}  If $A$ is
 quasidiagonal relative to $B$
 then there is a natural isomorphism
$$
Z_1(A,B) \,\cong\,   QD(A,B)   .
\tag 2.2
$$
\qed\endproclaim
\medbreak

Of course this implies that $QD(A,B)$ is a subgroup of $\KKgraded *.A.B $.

\medbreak

Here is the principal result of this paper.
\medbreak

\proclaim {Theorem 2.3} Suppose that $A \in \widetilde{\Cal N} $, $B$ is 
a separable $C^*$-algebra,  and
$A$ is quasidiagonal relative to $B$.
Then there is a
natural isomorphism  of topological groups
$$
   QD(A,B)  \,\cong\,  \pext  {K_*(A)}.{K_*(B)} _{0}
$$
  regarded as subgroups
of $\KKgraded 1.A.B   $ via
the canonical inclusion $\delta $ in the UCT.
\endproclaim
\medbreak

\medbreak\medbreak

\demo{Proof} We have isomorphisms of topological groups
$$
\align
QD(A,B) \,&\cong\, Z_1(A,B)
\intertext{by 2.2}
&\cong\,  \pext {K_*(A)}.{K_*(B)} _{0}
\endalign
$$
by Theorem 1.6. This completes the proof.\qed
\enddemo

 \medbreak

\proclaim {Corollary 2.4} Suppose that $A \in \widetilde{\Cal N} $  and
that $A$ is quasidiagonal. Then
$$
QD(A, \Cal K) \,\cong\, \pext {K_0 (A)}.{\Bbb Z}    .
$$
\qed\endproclaim
\medbreak

 Theorem 2.3 and Corollary 2.4 may be used readily since much is known
about computing the $Pext $ groups shown;    see \cite {22}.

\medbreak

\flushpar {\bf {Remark 2.5.\,\,}} The group $QD(A,B)$ may well have torsion.
For instance, suppose that $A \in \widetilde{\Cal N} $, $A$ is quasidiagonal
     and $K_*(A)$ is
torsionfree but not free. Then\footnote {since
$$
\pext G.H  \cong \ext G.H
$$
whenever $G$ is torsionfree.}
$$
QD(A,\Cal K) \,\cong\,
\pext {K_0(A)}.{\Bbb Z}       \,\cong\,
\ext {K_0(A)}.{\Bbb Z}         
$$
 which is always an uncountable, divisible group (cf. \cite{22, 
 Theorem 9.7}, due to  C.U. Jensen.)

Now choose $A \in \widetilde{\Cal N} $  with
$K_0(A) = \Bbb Z_p$, the integers localized at $p$.
Then
$$
QD(A,\Cal K) \,\cong\, \Bbb Q^{\aleph _0} \,\oplus\, \Bbb Z(p^{\infty }) .
$$
where  $\Bbb Z(p^{\infty }) $ denotes the $p$-torsion subgroup of $\Bbb Q/\Bbb Z$.
This point was overlooked in \cite {15, Cor. 5.4}. In that paper the
primary interest was the case with  $K_*(B)$   torsionfree and
$K_*(A)$   finitely generated.  If  $G$ is finitely generated then $\pext G.H   = 0$
for all $H$, so if $A \in \widetilde{\Cal N} 
 $ with $K_*(A)$ finitely generated then every
quasidiagonal extension  of the form
 $$
0 \to B\tensor\Cal K \to E \to A \to 0
$$
is trivial in $\KKgraded 1.A.B$   .

\medbreak

\proclaim{Corollary 2.6} Suppose that $A \in \widetilde{\Cal N} $,
 $A$ is quasidiagonal relative to $B$,
and either
\roster
\item  $K_*(A)$ is a direct sum of (finite and/or infinite)
cyclic groups.

or
\item $K_*(B)$ is algebraically compact\footnote{A group
is {\it algebraically compact} if it is a direct summand in
every group that contains it as a pure subgroup. Equivalently,
 it is algebraically compact if it is algebraically a direct
summand in a group which admits a compact topology. Examples
include compact groups, divisible groups, and bounded groups.
A group is algebraically compact if and only if it is of
the form
$$
D \oplus \prod _p D^p
$$
where $D$ is divisible and for each prime $p$ the group
$D^p$ is the completion in the $p$-adic topology
of the direct sum of cyclic $p$-groups and groups of
 $p$-adic integers. For any sequence of abelian groups
$\{G_i \}$ the group $\prod G_i /\oplus G_i $ is
algebraically compact and its structure is known.}.

\endroster
 Then any essential quasidiagonal extension
$$
0 \to B\tensor\Cal K  \to E \to A \to 0
$$
is a trivial extension.
\endproclaim

\medbreak
\demo{Proof} If $G$ is a direct sum of cyclic groups then
$\pext G.H = 0$ for all groups $H$.
Dually, if $H$ is algebraically compact then
$\pext G.H = 0$ for all groups $G$.
 The result is then immediate
from Theorem 2.3.
\qed\enddemo
\medbreak

\proclaim{Corollary 2.7} 
Suppose that $A \in \widetilde{\Cal N} $, $A$ is quasidiagonal relative 
to $B$, $K_*(A)$  is torsionfree, and we are given an
essential extension 
$$
\tau : \qquad   0 \to B\tensor\Cal K \to E \to A \to 0.
$$
 Suppose further that the connecting homomorphism
$$
K_*(A) \to K_{*-1}(B)
$$
is trivial. Then the extension is quasidiagonal.
If $K_*(A)$ is free then the extension is trivial.
\endproclaim
\medbreak
\demo{Proof} The fact that the connecting homomorphism
vanishes implies that the class of $\tau $ lies in
the group $\ext{K_*(A)}.{K_*(B)} _0 $. Since $K_*(A)$ is torsionfree
 we know \cite {22} that
$$
\usualext _0 \,=\, \usualpext _0
$$
and then we have
$$
[\tau ] \in   \usualpext _0 \,=\, QD(A,B)
$$
as required. If $K_*(A)$ is free then $\pext {K_*(A)}.H = 0$ for 
any abelian group $H$ and so $Q(A,B) = 0$ and the extension 
is trivial.\qed\enddemo

\flushpar {\bf Remark 2.8}. The group $\pext G.H $ is either
zero or huge. For example,
\medbreak
\roster
\item  (Warfield) If both $G$ and $H$ are
countable   groups
 with $\pext G.H  \neq 0$ and $H$ torsionfree, then
$\pext G.H $  
has $\Bbb Q^{\aleph _o}$ as a direct summand. It also
may have torsion.
\medbreak
\item (R. Baer) If $G$ is a torsionfree   group and $H$ is
a countable torsion group then
$$
\pext G.H = \Bbb Q^n
$$
with either $n = 0$ or $n \geq \aleph _o$.
\medbreak
\endroster

\medbreak
For proofs and references to these facts and for many more
examples, please see  \cite{22}.

\medbreak

\newpage\beginsection {3. Homomorphisms and Split Morphisms}

In this section the    identification of the quasidiagonal
elements of Theorem 2.3  is made concrete.
For instance, we  answer a  simple question. Suppose that
$f: A \to B$. Then when is the canonically 
associated class $\beta ^A[f] \in 
\KKgraded 1.{SA}.B $ 
a quasidiagonal class?
   There are certain easy algebraic
necessary conditions, and we show that (for $A \in \widetilde{\Cal N} $) these
conditions are sufficient. Finally, we demonstrate 
how to deal with an extension given explicitly.

In applications many of the most interesting classes in $KK$-theory
come from $*$-homo\-mor\-phisms $f: A \to B$. The class $[f] \in 
\KKgraded 0.A.B $ determines the structure of $\KKgraded *.A.B $
as a module over the ring $\KKgraded *.A.A $ and hence has 
special importance.

Let 
$$
\tau _{\Bbb C}\,\,:\qquad\qquad  
0 \to \Bbb C\tensor\Cal K \to \Cal T  \to  S\Bbb C  \to 0
$$
denote the universal Toeplitz extension which generates
the group 
$$
\KKgraded 1.{S\Bbb C}.{\Bbb C  } \cong \Bbb Z .
$$
We may tensor the extension with $B$ to obtain the extension
$$
\tau _B\,\,:\qquad\qquad
0 \to B\tensor\Cal K \to B\tensor\Cal T  \to  SB  \to 0.
$$
Note that 
$$
\gamma (\tau _B ) = 1 : K_*(B) \to K_*(B)
$$
and hence 
$$
\gamma (\tau _B ) = 0 \qquad\text{if and only if}\qquad  K_*(B) = 0,
$$ 
which would imply that $\KKgraded *.B.B = 0$ if $B \in \widetilde{\Cal N} $.
Thus if $B \in \widetilde{\Cal N} $ then
the extension $\tau _B $ is $B$-quasidiagonal if and only
if it is trivial, and this only happens when $K_*(B) = 0$. 

There is a commuting diagram
$$
\CD
\KKgraded 0.{\Bbb C}.{\Bbb C} @>{\beta ^{\Bbb C}}>> 
  \KKgraded 1.{S\Bbb C}.{\Bbb C} \\
@VV{\iota _B}V    @VV{\iota _B}V   \\
\KKgraded 0.{B}.{B} @>{\beta ^B}>> 
  \KKgraded 1.{SB}.{B}
 \endCD
 $$
 where $\beta ^A$ represents the Bott isomorphism 
 in the first variable of $KK$
 and $\iota _B$ the
 canonical structural map.  Then
 $$
 \align
 \tau _B &= \iota _B (\tau _{\Bbb C}) \\
 \intertext{by construction} 
 &= \iota _B\beta ^{\Bbb C}([1_{\Bbb C}]) \\
 \intertext{since this is the universal Toeplitz extension}  
 &=  \beta ^{B}\iota _B ([1_{\Bbb C}]) \\
 \intertext{since the diagram commutes} 
 &= \beta ^{B} ([1_{B}])  .
 \endalign
 $$

The
 diagram 
 $$
 \CD
 \KKgraded 0.B.B  @>{\beta ^B}>> \KKgraded 1.{SB}.B  \\
 @VV{f^*}V    @VV{(Sf)^*}V  \\
 \KKgraded 0.A.B   @>{\beta ^A}>> \KKgraded 1.{SA}.B  
 \endCD
 $$
 commutes by the naturality of the Bott isomorphism, and hence
 we conclude that
 $$
 \beta ^A([f]) = \beta ^Af^*([1_B]) = (Sf)^*\beta ^B([1_B]) 
  = (Sf)^*\tau _B
  \in  \KKgraded 1.{SA}.B  .
 $$

 Thus we may represent the class  $\beta ^A([f])$ as the pullback
 of the extension $\tau _B $ by $Sf$, namely
 $$
 \CD
  0 @>>>  B\tensor\Cal K  @>>> E  @>>>   SA  @>>>  0\\
  @. @VV{\cong}V    @VVV  @VV{Sf}V   @. \\ 
0 @>>>  B\tensor\Cal K  @>>> B\tensor\Cal T  @>>>   SB  @>>>  0 
\endCD
$$
(If $Sf$ is not mono then as usual we add on a trivial extension so that 
the resulting Busby map classifying the pullback is mono.)
So - here is a canonical extension associated to the map $f$. 
The $C^*$-algebra $SA$ is always quasidiagonal
and so the trivial extension added on to make $Sf$ mono was 
quasidiagonal. Thus it makes
sense to ask when the class $\beta ^A([f])$ is 
quasidiagonal.

 Note that 
 $$
 \gamma (\beta ^A [f] ) = (\beta ^A)^*\gamma ( f) = 
 (\beta ^A)^*f_*  
 $$
 where 
 $$
 (\beta ^A)^* : \hom{K_*(A)}.{K_*(B)}  \overset{\cong}\to\longrightarrow 
 \hom{K_*(SA)}.{K_*(B)}
 $$
and since $(\beta ^A)^* $ is an isomorphism (shifting parity, 
of course), we see that
$$
\gamma (\beta ^A [f] ) = 0 \qquad \text{if and only if} \qquad 
f_* = 0 : K_*(A) \to K_*(B).
$$

In order to state the answer, we must first introduce the 
mapping cone of $f$.
Recall that the  {\it { mapping cone}} $Cf$ is the $C^*$algebra
$$
Cf = \{ (\xi , a) \in B[0,1]\oplus A : \xi (0) = 0,\,\,\xi(1) = f(a) \}
$$
with associated
 mapping cone sequence (cf. \cite {17})
$$
0 \to  SB  \to Cf   \to  A   \to   0  .
$$
that has equivalence class 
$$
\beta _B [f] \in  \KKgraded 1.A.{SB}  
$$
where $\beta _B $ represents Bott periodicity in the second variable
of $KK$.
Note that the diagram 
$$
\CD
\KKgraded 0.A.B  @>{\beta ^A}>>   \KKgraded 1.A.{SB}  \\
 @AA{f_*}A    @AA{f_*}A   \\
 \KKgraded 0.A.A  @>{\beta ^A}>>  \KKgraded 1.{A}.{SA}
\endCD
$$
commutes and hence 
$$
\beta _B([f]) = \beta _B f_* ([1_A]) = (Sf)_*\beta _A  ([1_A])
$$
where $\beta _A  ([1_A])$ is the class of the mapping cone of 
the identity map $1 : A \to A$ which has the form
$$
0 \to SA \to CA  \to A  \to 0.
$$
This class is $KK$-invertible with $KK$-inverse the Toeplitz
class $\tau _A \in \KKgraded 1.{SA}.A $.\footnote{This fact is 
the core case in the proof of Bott periodicity in the $KK$-context.}

\proclaim {Proposition 3.1}  Suppose that $A \in \widetilde{\Cal N} $.
 Let $f: A \to B$ be a $*$-homomorphism.
Then the class $\beta ^A [f]  \in \KKgraded 1.{SA}.{B}  
  $ is quasidiagonal if and only if both of the following
conditions hold:
\roster
\item The induced homomorphism
$
f_* : K_*(A) \to K_*(B) $  is trivial; and,
\medbreak
\item
The associated short exact sequence
$$
0 \to  K_*(SB)  \to K_*(Cf)  \to  K_*(A)  \to  0
$$
is pure exact.
\endroster
If $K_*(A)$ is torsionfree then condition (2) is automatically satisfied, so
that $\beta ^A[f]$ is quasidiagonal if and only if $f_* = 0$.
\endproclaim
\medbreak

\medbreak
\demo {Proof} By virtue of Theorem 2.2 it suffices to determine
when the class $\beta ^A[f]$ is in $Z_1(SA,B)$. The Bott 
map is a homeomorphism, by (1.2), and hence it is equivalent to ask 
when the class $[f] \in \KKgraded 0.A.B $ lies in the 
subgroup $Z_0(A,B)$. 
 By Theorem 2.3  it suffices to show that
$$
[f] \in {\pext {K_*(A)}.{K_*(B)} _0} 
$$
 if and only if
conditions 1) and 2) hold. Condition 1) is necessary and sufficient
for
$$
[f] \in {\ext {K_*(A)}.{K_*(B)} } 
$$
by the UCT, and condition 2) picks out those elements of
${  \ext {    K_*(A)    }.{   K_*(B)    } }  $
which lie in ${\pext {K_*(A)}.{K_*(B)} }  $.
If $K_*(A)$ is torsionfree then
$$
 {\pext {K_*(A)}.{K_*(B)} }   \,\,= \,\, {\ext {K_*(A)}.{K_*(B)} }
$$
and the corollary follows.
\qed\enddemo
\

\medbreak

Suppose given an extension of $C^*$-algebras 
$$
0 \to B  \overset{i}\to\longrightarrow 
E \overset{p}\to\longrightarrow A  \to 0
$$
which is split by a $*$-homomorphism $s: A \to E$. 
Then there is a 
 unique class   $\pi _s \in \KKgraded 0.E.B $ called the 
 {\it{splitting morphism}}
 with
the property that
$$
i_*(\pi _s) = [1] - [sp]  \in \KKgraded 0.E.E       .
$$
Note that $\pi _s$ induces a homomorphism 
$$
\gamma(\pi _s ) : K_*(E) \longrightarrow K_*(B) 
$$
but this homomorphism generally does {\it{not}} arise from a $*$-homomorphism
$$
E \to B\tensor\Cal K . 
$$

Every $KK$-class may be represented
as  the $KK$-product of a class
induced by a homomorphism and by a splitting morphism 
\cite{1, \,\,\S\S 17.1.2, 17.8.3}.
So we wish to know when a splitting morphism corresponds to 
a quasidiagonal class. Here is the answer.

\medbreak
\proclaim {Proposition 3.2}  Suppose that $A \in \widetilde{\Cal N} $.
  Further, suppose given an extension
of $C^*$-algebras
$$
0 \to B  \overset{i}\to\longrightarrow E \overset{p}\to\longrightarrow A  \to 0
\tag {4.3}
$$
which is split by a $*$-homomorphism $s: A \to E$.  Let
$
\pi _s  \in \KKgraded 0.E.B  $ 
be the associated splitting morphism.
Consider the class 
$$
\beta ^E (\pi _s)
\in 
 \KKgraded 1.{SE}.{B}  .
$$
Then the following
are equivalent:
\roster
\item  $\beta ^E (\pi _s) $ is a quasidiagonal class.
\medbreak
\item     $K_*(B)  = 0$.
\medbreak
\item $\KKgraded *.A.B = 0$.
\medbreak
\item  $\pi _s  = 0$  .
\medbreak
\endroster
\endproclaim
\medbreak

\demo {Proof} The implications $3) \to  4) \to
1)$  are obvious.  We shall show
that $1) \to  2) \to  3)$.

Suppose first that $\beta ^E(\pi _s) $ is a quasidiagonal
class.
 Then $\gamma (\beta ^E(\pi _s)) = 0$
by Theorem 2.3, which implies that $\gamma  (\pi _s)  = 0$.
Then
 $$
 1 - s_*p_*  \,=\, \gamma ([1] - [sp] ) \,=\,
\gamma (i_*(\pi _s)) \,=\, i_*(\gamma (\pi _s)) \,=\, 0
$$
so that  $p_* : K_*(E) \to K_*(A) $ is an isomorphism and thus
$K_*(B) = 0$ as required. Thus $1) \to 2)$.
If $K_*(B) = 0$ then $\KKgraded *.E.B = 0$ by the UCT. Thus $2) \to 3)$.
\qed\enddemo

\medbreak
Proposition 3.1 tells when a class $\beta ^A[f]$ is quasidiagonal. 
Proposition 3.2
tells when 
the class $\beta ^A(\pi _s)$   is quasidiagonal. The following
theorem describes when a class with factorization
$$
x =  [f] \tensor _D \pi _s = f^*(\pi _s) 
$$
corresponds to a quasidiagonal class $\beta ^A(x)$.
 Every $KK$-class has this form.
Thus this theorem gives a complete solution to the quasidiagonality
of the associated class $\beta ^A(x) \in \KKgraded 1.{SA}.B $. 
We state the theorem for the class $x$ for simplicity, 
remembering that $\beta ^A(x)$ is a quasidiagonal class if 
and only if $x \in Z_0(A,B) \cong \usualpext _1$.

\proclaim {Theorem 3.3}  Suppose that $A \in \widetilde{\Cal N} $.
 Let $x \in \KKgraded 0.A.B $
with factorization $$
x =  [f] \tensor _D \pi _s = f^*(\pi _s)
$$  with respect
to the map
$$
f: A \to D
$$
 and   extension
$$
0 \to   B\tensor \Cal K \to    D   \overset{p}\to\longrightarrow A  \to 0
$$
with splitting  $s: A  \to D$ and splitting morphism
 $\pi _s  \in \KKgraded 0.D.B $.

\medbreak\flushpar Then:
 \medbreak
\flushpar (a)  The following conditions are equivalent:
\roster
\item   $$
       \gamma (x) = 0 \in \hom {K_*(A)}.{K_*(B)} _0   .
        $$  \medbreak
\item                                       $$
           \align
            Im(f_* : K_*(A) \to K_*(D) )
            &\subseteq Ker( \gamma (\pi _s)  : K_*(D) \to K_*(B)) \\
            &= Im (s_*: K_*(A) \to K_*(D) )   .
        \endalign
           $$
 \medbreak
\endroster
\flushpar (b)     Suppose that $\gamma (x) = 0$, so that
$$
x \in \ext {K_*(A)}.{K_*(B)}   _1 \subseteq \KKgraded 0.A.B  .
$$
  Then
$$
x = z \tensor _D \pi _s
\tag *
$$
for some $$
z \in \ext {K_*(A)}.{K_*(D)}   _1 \subseteq \KKgraded 0.A.D .
$$
The element $z$ is unique modulo the subgroup
$$
s_*(      \ext {   K_*(A)    }.{   K_*(A)   }     )_1.
$$
Conversely, any element $x$ of   form (*) is in the group
$\ext {K_*(A)}.{K_*(B)}   _1$.
\medbreak
\flushpar (c)  Suppose that
$x \in \pext {K_*(A)}.{K_*(B)}   _1 \subseteq \KKgraded 0.A.B  $.
 Then
$$
x = z \tensor _D \pi _s
\tag {**}
$$
for some
$$
z \in \pext {K_*(A)}.{K_*(D)}   _1 \subseteq \KKgraded 0.A.D .
$$
The element $z$ is unique modulo the subgroup
$$
s_*(      \pext {   K_*(A)    }.{   K_*(A)   }     )_1.
$$
Conversely, any element of form (**) is in the group
 $ \pext {   K_*(A)    }.{   K_*(B)   }   _1  $  .
\medbreak
\endproclaim

\demo {Proof} For Part (a) we compute:
$$
\gamma (x)   = \gamma (\, [f]  \tensor _D \pi _s \,)
       =  \gamma (\pi _s)f_*
 $$
 and hence
$\gamma (x) = 0 $ if and only if $Im(f_*) \subseteq Ker(\, \gamma (\pi _s)  \,)$. The
identification $Ker(\, \gamma (\pi _s)  \,) = Im (s_*)$ is immediate from the
definition of $ \gamma (\pi _s)  $. This proves Part (a).
\enddemo
\medbreak

In order to prove Parts (b) and (c)  the following Lemma is required.    It uses
the notation and assumptions of   Theorem 3.3.
\medbreak

\proclaim {Lemma 3.4} The diagram
$$
\CD
\ext {    K_*(A)    }.{     K_*(D) }       _\starone        @>{\delta _D}>>
     \KKgraded *.A.D   \\
@VV{ \gamma (\pi _s)  _* }V        @VV{ (-)\tensor _D  \pi _s }V    \\
Ext_{\Bbb Z}^1 (K_*(A),K_*(B) ) _\starone
          @>{\delta _B}>>             \KKgraded *.A.B
\endCD
$$
commutes, where the maps $\delta $ are the inclusion maps from the UCT.
\endproclaim
\demo {Proof} This does not follow immediately from the naturality of
the UCT, since the map $  (-)\tensor _D\pi _s $ is not induced by
 a map of $C^*$-algebras. We argue as follows. Expand the diagram to
the diagram
$$
\CD
0       @.    0    \\
@VVV      @VVV      \\
\ext {K_*(A)}.{K_*(A)}   _\starone @>{\delta _A}>>   \KKgraded *.A.A   \\
@VV{s_*}V    @VV{s_*}V    \\
Ext_{\Bbb Z}^1 (K_*(A),K_*(D) ) _\starone
     @>{\delta _D}>>               \KKgraded *.A.D   \\
@VV{ \gamma (\pi _s)  _* }V        @VV{ (-)\tensor _D  \pi _s }V    \\
\ext {K_*(A)}.{K_*(B)} _\starone   @>{\delta _B}>>   \KKgraded *.A.B    \\
@VVV   @VVV   \\
0        @.           0
\endCD
$$
\medbreak
Each column is split exact,
since $s$ is a splitting,
and the
 upper square commutes by the naturality of the UCT (since the map
$s_*$ {\it {is}} induced by the map $s: A \to D$) and it is easy to
see that the map $\delta _B $ is the quotient map, making the
lower square commute. \qed\enddemo

  Returning to the proof of Part (b), we see from   Lemma 3.4 that the
restriction of the map $(-)\tensor _D \pi _s$ to the
map
$$
(-)\tensor _D \pi _s :
\ext {K_*(A)}.{K_*(D)}  \to   \ext {K_*(A)}.{K_*(B)}
$$
is just the split surjection $ \gamma (\pi _s)_*  $. This implies that
 each element
of the topological group $\ext {K_*(A)}.{K_*(B)} _1$ is of the 
requisite form and the
indeterminacy is as stated.

For Part (c), consider the commuting diagram
\medbreak
$$
\CD
\pext {K_*(A)}.{K_*(D)}   @>{(-)\tensor _D \pi _s}>>  \pext {K_*(A)}.{K_*(B)} \\
@VV\psi V          @VV\psi V    \\
\ext {K_*(A)}.{K_*(D)}   @>{(-)\tensor _D \pi _s }>>   \ext {K_*(A)}.{K_*(B)}
\endCD
$$
 The lower horizontal map is just the map $  \gamma (\pi _s)  _* $ by the Lemma, and so the
 upper horizontal map (which is a further restriction, of course) is just the map
  $  \gamma (\pi _s)  _*  $  as well. Then the conclusion of the Theorem is an easy consequence
 of previous results and the fact that both horizontal maps are split
 surjections with known kernels. \qed
\medbreak

Here is our solution to the relative
quasidiagonality  problem, applied to a concrete extension.

\proclaim {Theorem 3.5}   Suppose that $A \in \widetilde{\Cal N} $ is quasidiagonal
relative to $B$. Suppose given an essential extension
$$
0 \to  B \tensor \Cal K   \to  E_\tau    \to A \to  0
$$
representing $\tau \in \KKgraded 1.A.B $
with corresponding representation 
$$
e_\tau : E_\tau \to \Cal L(H_B). 
$$
 Then
the representation $e_\tau $ is $B$-quasidiagonal
 if and only if both of the following conditions
hold:
\roster
\item $\gamma (\tau ) $ = 0, or, equivalently, the boundary homomorphism
$K_*(A) \to K_\starone (B)$ is trivial; and
\medbreak
\item
$$
\tau \in  \underset{n}\to\bigcap \, n\ext {K_*(A)}.{K_\starone (B)} _0
$$
or, equivalently,
$$
\tau \in Ker[\varphi : \ext {K_*(A)}.{K_*(B)} _0   \longrightarrow
 \ext {K_*(A)}.{K_*(B)} _0 {}^{\wedge}   ]
$$
where $G{}^{\wedge } $ denotes the $\Bbb Z$-adic completion of $G$.
\endroster
If $K_*(A)$ is torsionfree then condition (2) is satisfied
automatically, so that $e_\tau $ is a $B$-quasidiagonal
representation
if and only if the boundary homomorphism 
$K_*(A) \to K_\starone(B) $ is trivial.
\endproclaim
\medbreak
\demo {Proof} Conditions 1) and 2) are exactly the conditions that
guarantee that $\tau $ lies in the subgroup  $ QD(A,B) $.  \qed\enddemo


\newpage
\beginsection {4. Purity and Quasidiagonality}

 In this section we take advantage of standard results in infinite abelian
groups to deduce results on quasidiagonality.

\medbreak

\proclaim {Proposition 4.1}

\flushpar (a)  \quad Suppose that $H$ is a countable  abelian
group. Then the following are equivalent:
\roster
\item $\pext G.H  = 0 $ for all countable abelian groups $G$.
\medbreak
\item $\pext G.H  = 0 $ for the groups $G = \Bbb Q$  and  $\Bbb Q/\Bbb Z $.
\medbreak
\flushpar
\item $H$ is algebraically
 compact.
 \medbreak
\endroster
\flushpar
(b)\quad
   Suppose that $G$ is a countable abelian group.
 Then the following are equivalent:
\roster
\item  $\pext G.H  = 0 $ for all countable abelian groups $H$.
\medbreak
\item  $\pext G.H  = 0 $ for all countable direct sums of
cyclic groups $H$.
\medbreak
\item $G$ is the direct sum of cyclic groups.
\endroster
\endproclaim
\medbreak
\demo {Proof} First concentrate on Part (a). Of course (a1) implies
(a2).  The implication (a3) implies (a1) is immediate from \cite {11, 53.4}.
 The implication (a2) implies (a3) is the least obvious, but it is
also found in    \cite {11, page 232.}

Turning to part (b), the implication (b1) implies (b2) is trivial, and
the implication (b3) implies (b1) follows from \cite {11, 53.4} as well.
For the following argument that (b2) implies (b3) I am indebted
to John Irwin. Suppose that $G$ satisfies condition (b2). Let
 $\widetilde G$   be the free abelian group on  the (countable) set
$$
\{ [g] : g \in G  \}
$$
  modulo the relations given by
$$
n[g] = 0
$$
if $g \in G$ has order $n$.  There is an obvious surjection
$\widetilde G   \to G$ and hence a short exact sequence
$$
\Theta :\qquad\qquad
   0  \to  K  \to   \widetilde G   \to   G \to    0   .
$$
This sequence is pure, since every torsion element of $G$ lifts to a
torsion element of $\widetilde G$ of the same order. The group $\widetilde G$
is countable by construction, hence $K$ is countable.
Further, $\widetilde G$ is a direct sum of cyclic groups.
The group
$K$ is a subgroup of a direct sum of cyclic groups and
by Kulikov's theorem (cf.  \cite {11, 20.1} )  $K$ itself is a direct
sum of cyclic groups. Thus
$$
\Theta  \in \pext G.K   = 0
$$
and hence the extension $\Theta $  must be split. This implies that $G$ is
isomorphic to a subgroup of a direct sum of cyclic groups and hence (by
Kulikov's theorem)
is itself a direct sum of cyclic groups.
\qed\enddemo
\medbreak

 The following theorem is the $KK$-version of the preceding, purely
algebraic results.
First a bit of notation.
For $G$  any countable abelian group,  let $C_G$ be a separable commutative
$C^*$-algebra such that for each $j = 0,1$,
$$
K_j(C_G )   = G.
$$
Note that $C_G$ exists by geometric realization \cite{16}
and is unique up to $KK$-equiv\-alence  by the UCT.
Each $C_G$ is quasidiagonal, since it is commutative, and 
hence $C_G$ is $B$-quasidiagonal for all separable $B$.\footnote{If
 desired we may choose $C_G = C_G^0 \oplus C_G^1 $ where
$ K_i(C_G^j) = G $ if $i = j$ and   $ K_i(C_G^j) = 0 $ if $i \neq j$.
Then we could use the $C_G^j$ separately in parts (2) and (3) of 
Theorem 4.2.}  

\medbreak

\proclaim {Theorem 4.2}

\flushpar (a)\quad Suppose that $B$ is a separable $C^*$-algebra.
Then the following
are equivalent:
\roster
\item For each $A \in \widetilde{\Cal N} $ with $A$ quasidiagonal relative to $B$,
$$
QD(A,B) \,=\,  0 .
$$
\item  For  $G = \Bbb Q  $  and  $G  =  \Bbb Q / \Bbb Z  $,
$$
QD(C_G,B) \,=\,  0 .
$$
\medbreak
\item $K_*(B)$ is algebraically compact.
\endroster
\medbreak
\flushpar (b) Suppose given a quasidiagonal
 $C^*$-algebra $A \in \widetilde{\Cal N} $. Then the
following are equivalent:
\roster
\item For each separable $C^*$-algebra $B$,
$$
QD(A,B) \,=\,  0 .
$$
\medbreak
\item For $H$ any direct sum of cyclic groups,
$$
QD(A,C_H) \,=\,  0 .
$$
\medbreak
\item  $K_*(A)$ is the direct sum of cyclic groups.
\endroster
\endproclaim
\medbreak

\demo {Proof}  First consider (a). The implication (a1) implies (a2)
is immediate.  The implication (a3) implies (a1) is elementary, since
if $K_*(B)$ is algebraically compact then
$$
\pext G.{K_*(B)}  = 0
\tag *
$$
for all groups $G$  by 4.1(1).   If $A \in \widetilde{\Cal N} $  with
$A$ quasidiagonal relative to $B$ then
$$
QD(A,B) \cong \pext {K_*(A)}.{K_*(B)}
$$
by Theorem 2.3, and
$$
\pext {K_*(A)}.{K_*(B)}   \,=\,   0
$$
by (*), completing the argument.

Next we show that (a2) implies (a3).  The condition (a2) implies
that
$$
\pext G.{K_*(B)}  \,=\, 0
$$
for $G = \Bbb Q$ and  $ G =  \Bbb Q / \Bbb Z $ and then Theorem 4.1
implies that $K_*(B)$ is algebraically compact.
This completes the proof of part (a).

The proof of part (b) is quite similar, and we comment only
on the deep implication (b2) implies (b3). Condition (b2) together
with Theorem 1.2  imply
that
$$
\pext {K_*(A)}.H   \,=\,0
$$
whenever $H$ is a direct sum of cyclic groups, and then Theorem
4.1(1) implies that $K_*(A)$ is itself a direct sum of cyclic groups.
This completes the proof of Theorem 4.2.
\qed\enddemo
\medbreak
\newpage


\beginsection {5. A problem of L.G. Brown }

Let $\theta : K_*(A)_t \to K_*(A) $
be the canonical inclusion of the torsion subgroup of $K_*(A)$.
L.G. Brown \cite {4,  page 63} showed that (with $B = \Cal K$) 
there is a relation between
quasidiagonality and the   kernel of the induced map
$$
\theta ^*  : \ext {K_*(A)}.{K_*(B)} \longrightarrow \ext { K_*(A)_t}.{K_*(B)}  .
$$
We generalize his result as follows.
Let $G_t$ denote the torsion subgroup of a group $G$ and $G_f = G/G_t $
denote the  maximal torsionfree quotient of $G$.
Recall \cite {23} that given $A$, there is an associated extension
 of
$C^*$-algebras
$$
0 \to A\tensor \Cal K  \to A_f \to   SA_t \to 0
$$
whose $K$-theory long exact sequence degenerates to the
pure short exact sequence
$$
0 \to
 K_*(A)_t \overset\theta\to\longrightarrow  K_*(A)  \to K_*(A)_f  \to 0  .
\tag *
$$
In particular,
$$
K_*(A_t) \cong K_*(A)_t  \qquad   {\text{and}} 
\qquad  K_*(A_f) \cong K_*(A)_f  .
 $$
Further, if $A \in \widetilde{\Cal N} $ then so are both 
$A_t$ and $A_f$.\footnote{We established this in  \cite{23}, Theorem 1.1, 
for $A \in \Cal N$ but once again it is clear by inspection of 
that proof that the statment holds for $A \in \widetilde{\Cal N} $.}
\medbreak
\medbreak

\proclaim {Theorem 5.1} Suppose that $A \in \widetilde{\Cal N} $      and $A$ is
quasidiagonal relative to $B$. Then:
\roster
\item   There is a natural commutative diagram with exact columns:
\medbreak
$$
\CD
0   @.    0   \\
@VVV   @VVV   \\
 Ker (Q\theta ^*)         @>\cong >>   Ker (\theta ^*)  \\
@VVV    @VVV   \\
 QD(A,B)    @>\cong >>    \pext {K_*(A)}.{K_*(B)} _0 \\
@VV{Q\theta ^*}V    @VV{\theta ^*}V   \\
 QD(A_t,B)    @>\cong >>   \pext {K_*(A)_t }.{K_*(B)} _0  \\
@VVV    @VVV   \\
0  @.  0
\endCD
$$
\medbreak
\item There is a natural   exact sequence
$$
  \hom {K_*(A)_t }.{K_*(B)} _0
\overset{\delta '}\to\longrightarrow
 QD(A_f,B) \to  Ker (\theta ^*)   \to 0
$$
where $\delta '$ is the boundary map in the $Hom$-$Pext$ long exact
sequence associated to the pure short exact sequence (*).
\medbreak
\item If  $\Cal Im(\delta ') = 0 $\footnote {This condition usually
holds. For instance, it holds if  $K_*(A)_t  $ is a direct
summand of $K_*(A)$ or if  $K_*(B)$ is torsionfree, and of
course it holds if $\hom {K_*(A)_t}.{K_*(B)} = 0$. }
then there is a natural isomorphism
$$
QD(A_f, B) \,\cong\,   Ker (\theta ^*)  .
$$
\medbreak\endroster\endproclaim

\medbreak
\demo {Proof} The short exact sequence  (*) is pure exact and
thus produces a six term $Hom$-$Pext$ sequence in the first variable.
Identifying entries using Theorem 2.3, one obtains the left 
column below, which is exact. The six term $Hom$-$Ext$ exact
sequence contributes the right column below, and this column
is also exact. 
$$
\CD
\hom {K_*(A)_t }.{K_*(B)} _0    @>\cong >>    \hom {K_*(A)_t }.{K_*(B)} _0   \\
@VV{\delta '}V    @VV{\delta }V    \\
QD(A_f,B)     @>\cong >>     \ext  {K_*(A)_f }.{K_*(B)} _0  \\
@VVV    @VVV    \\
QD(A,B)     @>>>     \ext  {K_*(A)}.{K_*(B)} _0  \\
@VV{Q\theta ^*}V    @VV{\theta ^*}V    \\
QD(A_t,B)     @>>>     \ext  {K_*(A)_t }.{K_*(B)} _0  \\
@VVV    @VVV    \\
0    @.   0
\endCD
$$
\medbreak\flushpar
The diagram commutes,
the horizontal maps are injections, 
 and 
 the map
$$
QD(A_f,B)    \to    \ext  {K_*(A_f)}.{K_*(B)} _0
$$
is an isomorphism since  $  K_*(A_f) $  is torsionfree.
 An easy diagram chase shows that $Ker(Q\theta ^*) \cong Ker(\theta ^*) $
from which 1) is immediate. Part 2) follows from expressing $Ker(Q\theta ^*) $
as the quotient of $QD(A_f,B) $ modulo the group $\Cal Im (\delta ' ).$
\qed\enddemo

\medbreak


\newpage
\beginsection {6.  Quasidiagonality and Torsion}

One of the   early applications of the Brown-Douglas-Fillmore theory was
contained in work   of L.G. Brown \cite {4}. He exhibited an example
of a bounded operator $T$ which was not quasidiagonal but such that
$T \oplus T$ was quasidiagonal. In fact Brown showed that $T\oplus T$
generated a trivial extension.

Here is an analysis of such behavior from our perspective.

\medbreak
\proclaim {Theorem 6.1} Suppose that $A \in \widetilde{\Cal N} $, $A$ is 
quasidiagonal relative to $B$,
 and $x \in \KKgraded 1.A.B $.
Then:
\roster
\item If $\gamma (x) \neq 0 $ then $x$ is not a quasidiagonal class.

\medbreak\item  If  
$\gamma (x)$ has
infinite order 
in the group $\usualhom $
then no multiple of $x$ is a quasidiagonal class.

\medbreak\item If $\gamma (x) \neq 0 $   and  $\gamma (x)$ has order
$n$ 
in the group $\usualhom $
then $knx$ is  a quasidiagonal class for each $k \in \Bbb N$.

\medbreak
\item If $\gamma (x) = 0$ then $x$ is a quasidiagonal
class  if and only if it is in
the kernel of the  natural  map
$$
\varphi : \ext {K_*(A)}.{K_*(B)} \longrightarrow 
  \ext {K_*(A)}.{K_*(B)}  \sphat   .
$$
\medbreak

\item If $K_*(A)$ is torsionfree then $x$ is 
a quasidiagonal class if and only
if $\gamma (x) = 0$.
\endroster
\endproclaim

\medbreak

\demo {Proof} This is all immediate from Theorem 2.3. \qed\enddemo

\medbreak

We apply Theorem 6.1 to the setting of $x \in \KKgraded 1.A.{\Bbb C} $.

\proclaim{Theorem 6.2} Suppose that $A \in \widetilde{\Cal N} $ 
and $A$ is a quasidiagonal $C^*$-algebra.
 Suppose given an essential extension
$$
\tau : \qquad         0 \to \Cal K \longrightarrow E_\tau  
\overset{p}\to\longrightarrow A \to 0
$$
representing
 $\tau  \in \KKgraded 1.A.{\Bbb C} $ so that
 $$
 \gamma (\tau ) : K_1(A) \longrightarrow K_0(\Cal K ) \cong  \Bbb Z.
 $$ 
 Then:
\medbreak
\roster
\item  If $\gamma (\tau ) \neq 0$ then $\tau $ is not a
quasidiagonal extension
and no multiple of $\tau $ is a quasidiagonal extension.
\medbreak
\item If $\gamma (\tau ) =  0$ then $\tau $ is a quasidiagonal extension
if and only if the short exact sequence
$$
0 \to \Bbb Z \longrightarrow K_0(E_\tau ) 
\overset{p_*}\to\longrightarrow K_0(A) \to 0
$$
is a {\it {pure}} short exact sequence.
\medbreak
\item If $K_0(A) $ is torsionfree then $\tau $ is quasidiagonal 
if and only if $\gamma (\tau ) = 0$.
\medbreak
\item If $K_0(A) $ is a 
direct sum of (finite and/or infinite) cyclic groups
 then 
$\tau $ is a quasidiagonal extension if and only if it is a trivial
extension.
\endroster
\endproclaim

\medbreak

\demo{Proof} The group $\hom {K_1(A)}.{\Bbb Z}$  is torsionfree
since $\Bbb Z$ is torsionfree. Thus 6.1(1 and 2) imply 1). 
Part 2) is a restatement of 6.1(4). Part 3) follows from 6.1(5). 
Part 4) holds since $\pext G.H = 0$ whenever $G$ is a direct
sum of cyclic groups.\qed
\enddemo

\flushpar {\bf {Remark 6.3.}} There are examples \cite {19 } where
the surjection
$$
\varphi : \ext {K_*(A)}.{K_*(B)} \longrightarrow
  \ext {K_*(A)}.{K_*(B)}  \sphat
$$
 is not a
{\it{split}} surjection. In fact,  using the example of
Christensen-Strickland cited there, it is possible to produce
 an extension
 $\tau  \in \KKgraded 1.A.B  $ which
satisfies all of the following conditions:
\roster
\item  $\gamma (\tau ) = 0$,
\medbreak
\item $\varphi (\tau ) \neq 0 $ so that $\tau $
 is  not a quasidiagonal class, and
\medbreak
\item For some $k \in \Bbb N$,
 $k\tau $ is a quasidiagonal class but  is not trivial.
\endroster
This sort of phenomenon can occur only when both $K_*(A)$ and
$K_*(B)$ have $p$-torsion for some fixed prime $p$ and then only
rarely.
\medbreak
\medbreak

Theorem 6.1 implies an early result \cite {6} of Brown, Douglas,
 and Fillmore 
  (BDF)
on quasidiagonality. Recall that a bounded operator $T$
 on a separable Hilbert space $H$ is 
{\it{essentially normal}} if $T^*T - TT^*$ is compact. The 
{\it{essential spectrum}} of $T$ is the spectrum of $\pi T \in \Cal Q(H)$.
We let $ind(T)$ denote the Fredholm index of the operator $T$.

\proclaim {Theorem 6.4 (BDF)} Suppose 
that $T \in \Cal L(H)$ is 
an essentially normal operator.
Let $X$ denote the essential spectrum of $T$.
Then the
following are equivalent:
\roster
\item $T$ is a quasidiagonal operator.
\medbreak
\item For each $\lambda \in \Bbb C - X$,
$$
ind (T - \lambda I ) = 0.
$$
\medbreak\item $T$ is of the form (normal) + (compact).
\endroster
\endproclaim
\medbreak

\demo{Proof}  There is a natural extension
$$
0 \longrightarrow \Cal K \longrightarrow 
C^*\{T, \Cal K \}  \longrightarrow C(X) \longrightarrow  0
$$
which gives rise to an element
$$
[T] \in  \KKgraded 1.{C(X)}.{\Bbb C}   .
$$
The operator $T$ is quasidiagonal if and only if $[T]$ is 
a quasidiagonal class, by 2.2 and the remark preceding it.
The fact that  $X \subset \Bbb C$
implies that $K_*(C(X) \cong K^*(X)$ is torsionfree. Hence $[T]$
is  a quasidiagonal class 
if and only if $\gamma ([T]) = 0$ by Theorem 6.1. So it suffices
to compute $\gamma ([T])$.
 BDF \cite{6}
   show that there is a natural isomorphism
$$
  \KKgraded 1.{C(X)}.{\Bbb C}
  \overset\gamma\to\longrightarrow \hom {K^1(X)}.{\Bbb Z}
\,\cong\, \widetilde H^0(\Bbb C - X )
$$
and this map takes $[T]$ to the function
$$
\lambda \longmapsto      ind (T - \lambda I )
$$
so that $[T] = 0$ if and only if $ ind (T - \lambda I ) = 0$ where defined.
This completes the proof.\qed\enddemo


\newpage

\beginsection {7. Lifting Quasidiagonality }

In this section we present a converse to the following theorem of Davidson,
Herrero, and Salinas.
\medbreak

\proclaim {Theorem 7.1} \cite {10}  Suppose
that
$$
\tau: \qquad\qquad 0 \longrightarrow \Cal K  \longrightarrow E_\tau  
 \longrightarrow  A   \longrightarrow  0
\tag {7.2}
$$
is an essential extension with associated faithful representation
$$
e_\tau : E_\tau \longrightarrow \Cal L(H) 
$$
 and suppose further that $E_\tau $ is separable 
 and nuclear. If $e_\tau $ is a quasidiagonal 
  representation
then   $A$ is a quasidiagonal $C^*$-algebra.
\endproclaim

There is an obvious obstruction to a converse to this theorem,
 known already to Halmos \cite {12}.
Suppose that $S$ is the unilateral shift.   Then   there is
a canonical associated  extension
$$
\tau _S : \qquad\qquad  0 \longrightarrow \Cal K   \longrightarrow
  C^*\{S,\Cal K \} \longrightarrow  A  \longrightarrow 0
$$
and the quotient   $A  \cong C(S^1)$ is commutative, hence
quasidiagonal. However the $C^*$-algebra
$C^*\{S,\Cal K \}$ itself is not quasidiagonal 
since it contains the unilateral
shift $S$, which
 is not quasidiagonal. Halmos demonstrates this by observing that
$S$ has non-trivial Fredholm index,
 whereas any Fredholm quasidiagonal operator must
have trivial Fredholm index.  In modern  jargon, the index map
$$
\gamma  : \KKgraded 1.{C(S^1)}.{\Cal K}  \longrightarrow
 \hom {K^1(S^1)}.{\Bbb Z} \,\cong\, \Bbb Z
$$
satisfies
$$
\gamma (\tau _S)(z) = -1  \neq 0
$$
 which is an obstruction to quasidiagonality.

Here is the complete story, at least within the 
  category $\widetilde{\Cal N}$.
\medbreak

\proclaim {Theorem 7.3} Suppose given the essential
  extension (7.2) with
$A   \in \widetilde{\Cal N} $
and quasidiagonal.
  Then
the representation $e_\tau : E_\tau \to \Cal L(H) $ 
 is a quasidiagonal representation  
 if and only if the following two conditions hold:
\medbreak
\roster
\item $\gamma (\tau ) = 0$, and
\medbreak
\item The resulting $K$-theory short exact sequence
$$
0 \to \Bbb Z  \to  K_0(E_\tau )   \to   K_0(A )  \to  0
$$
is a pure exact sequence.
\endroster
If in addition $K_0(A  ) $ is torsionfree then condition 2) is automatically
satisfied, so that $e_\tau $ is a quasidiagonal representation
 if and only if $\gamma (\tau ) = 0$.
\endproclaim
\medbreak
\demo {Proof} Theorem 2.3 reduces in this case to the identification
$$
QD(A  , \Cal K) \cong  \pext {K_0(A  )}.{\Bbb Z}  .
$$
Condition 1) is equivalent to
$$
\tau \in \ext {K_0(A  )}.{\Bbb Z}
$$
and condition 2) is simply a statement that
$$
\tau \in \pext {K_0(A  )}.{\Bbb Z}   .
$$
If in addition $K_0(A  )$ is torsionfree then
$$
 \pext {K_0(A  )}.{\Bbb Z}  \,\cong\,  \ext {K_0(A  )}.{\Bbb Z}   .
$$
This completes the proof.\qed\enddemo
\medbreak

\flushpar{\bf {Remark 7.4.}}  It is instructive to compare 
our Theorem 7.3 with a related result due to N. Brown
and M. Dadarlat \cite{5}
where they overlap.
 Consider the essential extension
$$
\tau: \qquad\qquad 0 \longrightarrow \Cal K  \longrightarrow
 E_\tau   \longrightarrow  A   \longrightarrow  0
$$
with associated faithful representation
$$
e_\tau : E_\tau \longrightarrow \Cal L(H) 
$$
and suppose that $A \in \widetilde{\Cal N} $ with $A$  quasidiagonal. 
Brown and Dadarlat show (\cite{5, Theorem 3.4}
 specialized to this case) that if $\gamma (\tau ) = 0$
 then $E_\tau $ is also quasidiagonal. We show that
 if $\gamma (\tau ) = 0$ {\it {and}} the resulting 
 $K$-theory short exact sequence is {\it{pure}} then the 
 representation $e_\tau $
 is a quasidiagonal representation. 
If $\gamma (\tau ) = 0 $
but
 the $K$-theory short exact sequence is {\it{not}} pure
then we conclude that even though $E_\tau $ {\it{is}} quasidiagonal,
the representation $e_\tau $ is {\it{not}} a 
quasidiagonal representation. 

Here is an example.  Let $G$ be any countable
torsion group. Then 
$$
\ext G.{\Bbb Z} \,\cong\, \hom G.{{\Bbb Q}/{\Bbb Z}} = \Cal P(G)
$$
the Pontrjagin dual group of $G$.\footnote{ For example,
if $G$ is finite then $\Cal P(G) = G$. If $G = \oplus_1 ^\infty \Bbb Z /p $,
the sum of countably many copies of the group $\Bbb Z/p $, then
$\Cal P(G) = \prod _1 ^ \infty \Bbb Z /p $ which is, of course,
uncountable. Note that $\Cal P(\Cal P(G)) = G$ by the Pontrjagin
 duality theorem, and hence if $G \neq 0$ then
  $ \ext G.{\Bbb Z}  \neq 0$. }

Choose a commutative $C^*$-algebra  
$A$ with $K_0(A) = G$ and $K_1(A) = 0$. 
This is always possible, and $A$ is unique up to 
$KK$-equivalence. Then 
$$
\hom {K_*(A)}.{\Bbb Z} = 0
$$
and hence the index map $\gamma $ is identically zero. 
Thus there is a natural isomorphism
$$
\KKgraded 1.A.{\Cal  K}  \,\cong\, \ext{K_0(A)}.{\Bbb Z} \,\cong\,
\ext G.{\Bbb Z}  \,\cong \, \Cal P(G).
$$
Using Brown and Dadarlat's result, we conclude that if 
$$
0 \to \Cal K  \to E_\tau   \to A   \to 0
$$
is any essential extension then the $C^*$-algebra $E_\tau $ 
is quasidiagonal. On the other hand, 
$$
QD(A,\Cal K) \cong \pext {K_0(A)}.{\Bbb Z} = 
\pext {G}.{\Bbb Z} = 0
$$
since $G$ is a torsion group and $\Bbb Z$ is torsionfree,
by \cite{22, Theorem 9.1}. Thus among all of the various
$\tau \in \KKgraded 1.A.{\Cal K} \cong \Cal P(G)$
  and associated 
representations 
$$
e_\tau : E_\tau \longrightarrow \Cal L(H) ,
$$
the only representation  $e_\tau $  that is quasidiagonal is the one 
corresponding to the trivial extension, where 
$$
[\tau ] = 0 \in \KKgraded 1.A.{\Cal K} .
$$

This also illustrates the phenomenon discovered by L.G. Brown
and discussed in Section 6, since any non-trivial extension
$\tau \in \KKgraded 1.A.{\Cal K} $ will have the property that
it itself is not quasidiagonal, but when added to itself 
enough times it becomes quasidiagonal and trivial.

\vglue 2in
\flushpar
  fiii071001preprint
 
 \end

\end